\documentclass[12pt]{article}
\input isolatin1.sty
\usepackage[latin1]{inputenc}
\usepackage[english]{babel}
\usepackage{amssymb,amsmath}
\usepackage{graphicx}
\newtheorem{cbthm}{Theorem}[section]
\newtheorem{cbcor}[cbthm]{Corollary}
\newtheorem{cblem}[cbthm]{Lemma}

\newtheorem{cbrem}[cbthm]{Remark}
\renewcommand{\a}{\alpha}

\title{Fibonacci numbers and orthogonal polynomials}

\author{Christian Berg}
 
\begin{document}
\maketitle

\begin{abstract} We prove that the sequence $(1/F_{n+2})_{n\ge 0}$ of
reciprocals of the Fibonacci numbers is a moment sequence of a
 certain discrete  probability, and we identify the
  orthogonal polynomials as little $q$-Jacobi polynomials with 
$q=(1-\sqrt{5})/(1+\sqrt{5})$. We prove that the corresponding kernel
polynomials have integer coefficients, and from this we deduce 
that the inverse of the corresponding Hankel matrices $(1/F_{i+j+2})$
have integer entries. We prove analogous results for the Hilbert matrices.
\end{abstract}



\section{Introduction}

In \cite{Ri} Richardson noticed that the {\it Filbert
matrices}
\begin{equation}\label{eq:filbert}
\mathcal F_n=\left(1/F_{i+j+1}\right),\quad 0\le i,j\le n,\quad n=0,1,\ldots,
\end{equation}
where $F_n,n\ge 0$ is the sequence of Fibonacci numbers, have the
property that all elements of the inverse matrices are integers. The
corresponding property for the {\it Hilbert matrices} $(1/(i+j+1))$ has
been known for a long time, see \cite{Co},\cite{Ch},\cite{S:L}. The
last reference contains a table of the inverse Hilbert matrices up to
$n=9$. 

Richardson gave 
an explicit formula for the elements of the inverse Filbert matrices and proved
it using computer algebra. The formula
shows a remarkable analogy with the corresponding formula for the
elements of the inverse Hilbert matrices in the sense that one shall 
replace some
binomial coefficients $\binom{n}{k}$ by the analogous {\it Fibonomial
 coefficients}
\begin{equation}\label{eq:fibonomial}
\binom{n}{k}_{\mathbb F}=\prod_{i=1}^k\frac{F_{n-i+1}}{F_i},\quad 0\le
k\le n,
\end{equation}
with the usual convention that  empty products are defined as 1.
These coefficients are defined and studied in \cite{Kn} and are
integers.  The sequence
 of Fibonacci numbers is $F_0=0, F_1=1,\ldots,$ with the recursion formula
$F_{n+1}=F_n+F_{n-1},\;n\ge 1$.

The Hilbert matrices are the Hankel matrices $(s_{i+j})$ corresponding to the
moment sequence 
$$
s_n=1/(n+1)=\int_0^1 x^n\,dx,
$$
and that the reciprocal matrices have integer entries can easily be
explained by the fact the corresponding orthogonal polynomials, namely
the Legendre polynomials, have integer coefficients. See section 4 for
details.

The purpose of the present paper is to show that 
$(1/F_{n+2})_{n\ge 0}$ is the moment sequence of a certain discrete
probability. Although this is a simple
consequence of Binet's formula for $F_n$, it does not seem to have
been noticed in the literature, cf. \cite{Ko}. We find the
corresponding
probability  measure to be  
\begin{equation}\label{eq:fibmea}
\mu=(1-q^2)\sum_{k=0}^\infty q^{2k}\delta_{q^k/\phi},
\end{equation}
where we use the notation
\begin{equation}\label{eq:golden}
\phi=\frac{1+\sqrt{5}}{2},\quad
q=\frac{1-\sqrt{5}}{1+\sqrt{5}}=\frac{1}{\phi}-1,
\end{equation} 
and $\delta_a$ denotes the probability measure with mass 1 at the point $a$.
The number $\phi$ is called the golden ratio.

The corresponding orthogonal polynomials are little $q$-Jacobi
polynomials
\begin{equation}\label{eq:lqJacobi}
p_n(x;a,b;q)={}_2\phi_1\left(\begin{matrix}q^{-n},abq^{n+1}\\aq\end{matrix};q,xq\right),
\end{equation}
see \cite{G:R},  specialized to the parameters $a=q,b=1$, with $q$
taking the value from (\ref{eq:golden}).

To be precise we define
\begin{equation}\label{eq:fibpol}
p_n(x):=F_{n+1}p_n(x\phi;q,1;q),
\end{equation}
and these polynomials have integer coefficients, since
they can be written
\begin{equation}\label{eq:fibpol1}
p_n(x)=\sum_{k=0}^n (-1)^{kn-\binom{k}{2}}\tbinom{n}{k}_{\mathbb F}
\tbinom{n+k+1}{n}_{\mathbb F}x^k.
\end{equation}
The orthonormal polynomials with respect to $\mu$ and having positive
 leading coefficients are given as
\begin{equation}\label{eq:fibonp}
P_n(x)=(-1)^{\binom{n+1}{2}}\sqrt{F_{2n+2}}p_n(x),
\end{equation}
so the kernel polynomial
$$
K_n(x,y)=\sum_{k=0}^n P_k(x)P_k(y),
$$
is a polynomial in $x,y$ with integer coefficients. If we denote 
$a_{i,j}^{(n)}$ the coefficient to $x^{i}y^{j}$ 
in the kernel polynomial, then it is  a general fact that the matrix
\begin{equation}\label{eq:matrixAn}
A_n=(a^{(n)}_{i,j}),\quad 0\le i,j\le n
\end{equation}
is the inverse of the Hankel matrix of the
problem $(s_{i+j})_0^n$, see Theorem \ref{cbthm:A} below.

This explains that the elements of the inverse of the matrix
 $(1/F_{i+j+2})_0^n$ are integers, 
 and we derive a formula for the entries from the
 orthogonal polynomials.

 The Filbert matrices (\ref{eq:filbert}) are not
 positive definite but non-singular, and they are the Hankel matrices of the
 moments of a (real-valued) signed measure with total mass 1. The orthogonal
 polynomials for this signed measure are the little $q$-Jacobi polynomials
\begin{equation}\label{eq:lqJacobi1}
p_n(x\phi;1,1;q)=\sum_{k=0}^n(-1)^{kn-\binom{k}{2}}
 \tbinom{n}{k}_{\mathbb F}\tbinom{n+k}{n}_{\mathbb F}x^k,
\end{equation}
and a simple modification of the positive definite case leads to
Richardson's  formula for the entries of the inverse of the 
Filbert matrices.  

The two results can be unified in the statement that for each
$\a\in\mathbb N=\{1,2,\ldots\}$ the sequence $(F_{\a}/F_{\a+n})_{n\ge 0}$ is a moment
sequence of a real-valued measure $\mu_{\a}$ with total mass 1. It is a
positive measure when $\a$ is even, but a signed measure when $\a$ is
odd. The orthogonal polynomials are little
$q$-Jacobi polynomials $p_n(x\phi;q^{\a-1},1;q)$. This is proved  in
section 3.

In section 2 we recall some basic things about orthogonal polynomials
both in the positive definite and in the quasi-definite case, and
Theorem \ref{cbthm:A} about the inverse of the Hankel matrices is proved.

In section 4 we briefly discuss the matrices $(1/(\a+i+j))_0^n$, where
$\a>0$. They are related to Jacobi polynomials 
transfered to the interval $]0,1[$ and belonging to the 
parameters $(0,\a-1)$. This leads to the formula (\ref{eq:choi}),
which for $\a=1$ is the formula for the elements of the inverse
Hilbert matrices.

After the circulation of a preliminary version of this paper (dated
April 10, 2006), Ismail
has extended the results of section 3 to a one parameter
generalization of the Fibonacci numbers, cf.  \cite{Is1}.

\section{Orthogonal Polynomials}

We start by recalling some simple facts from the theory of orthogonal
polynomials, cf. \cite{Ak} or \cite{Is} and in particular \cite{Chi}
for the quasi-definite case.

\medskip
{\it The positive definite case}.

We consider the set $\mathcal M^*$ of probability measures on $\mathbb
R$ with moments of any order and with infinite support. The moment
sequence of $\mu\in\mathcal M^*$ is 
\begin{equation}\label{eq:mom}
s_n=s_n(\mu)=\int x^n\,d\mu(x),\quad n=0,1,\ldots,
\end{equation}
and the corresponding Hankel matrices are given by
\begin{equation}\label{eq:Hankel}
H_n=\begin{pmatrix} s_0 & s_1 & \cdots & s_n\\
s_1 & s_2 & \cdots & s_{n+1}\\
\vdots & \vdots & & \vdots\\
s_n & s_{n+1} & \cdots & s_{2n}
\end{pmatrix},\quad n=0,1,\ldots.
\end{equation}
The orthonormal polynomials $(P_n)$ for $\mu$ are uniquely determined
by the equations
\begin{equation}\label{eq:orthpol}
\int P_n(x)P_m(x)\,d\mu(x)=\delta_{n,m},\quad n,m\ge 0,
\end{equation} 
and the requirement that $P_n$ is a polynomial of degree $n$ with
positive leading coefficient. This coefficient is equal to
\begin{equation}\label{eq:leading}
\sqrt{D_{n-1}/D_n},
\end{equation}
where $D_n=\det H_n$.
The reproducing kernel for the polynomials of degree $\le n$ is
defined as
\begin{equation}\label{eq:kernel}
K_n(x,y)=\sum_{k=0}^n P_k(x)P_k(y),
\end{equation} 
and is called the kernel polynomial. It is clear that we can write
\begin{equation}\label{eq:A}
K_n(x,y)=\sum_{i=0}^n\sum_{j=0}^n a^{(n)}_{i,j}x^{i}y^{j},
\end{equation}
where the numbers $a^{(n)}_{i,j}$ are uniquely determined and satisfy
$a^{(n)}_{i,j}=a^{(n)}_{j,i}$. If we collect these numbers in 
an $(n+1)\times(n+1)$-matrix
$A_n=(a^{(n)}_{i,j})$, then it is the inverse of the  Hankel matrix
$H_n$:

\begin{cbthm}\label{cbthm:A} 
$$
A_nH_n=H_nA_n=I_n,
$$
where $I_n$ is the unit matrix of order $n+1$.
\end{cbthm}

{\it Proof}. For $0\le k\le n$ we have by the reproducing property
\begin{equation}\label{eq:AH}
\int x^kK_n(x,y)\,d\mu(x)=y^k.
\end{equation}
On the other hand we have
$$
\int x^kK_n(x,y)\,d\mu(x)=
\sum_{j=0}^n(\sum_{i=0}^n s_{k+i}a^{(n)}_{i,j})y^j,
$$
and therefore
$$
\sum_{i=0}^n s_{k+i}a^{(n)}_{i,j}=\delta_{k,j}.
$$
 \quad$\square$

\medskip
{\it The quasi-definite case}.

If $\mu$ is a real-valued signed measure on $\mathbb R$ with total
mass 1 and moments of any order, one can still define the moments
 (\ref{eq:mom}) and the corresponding Hankel matrices
 (\ref{eq:Hankel}). To define orthogonal polynomials one has to assume
 that (\ref{eq:Hankel}) is a non-singular matrix for any $n$,
 i.e. that the determinants satisfy $D_n=\det H_n\ne 0$.
On the other hand, if orthogonal polynomials exist with respect to a
signed measure, then the Hankel determinants are non-zero.
 See \cite[Theorem 3.1]{Chi} for
 details. In this case the orthonormal polynomial $P_n$ is uniquely
determined by the requirement that the leading coefficient
$\sqrt{D_{n-1}/D_n}$ is either positive or purely imaginary with
positive imaginary part.
The corresponding kernel polynomial $K_n$ has real coefficients, and
 Theorem \ref{cbthm:A} remains valid.

\section{Fibonacci numbers}

The Fibonacci numbers can be given by the formula
\begin{equation}\label{eq:Binet}
F_n=\frac{1}{\sqrt{5}}(\phi^n-{\hat\phi}^n),\quad n\ge 0
\end{equation}
usually called Binet's formula, but it is actually older, see
\cite{Kn},\cite{Ko}. Here 
$$
\phi=\frac{1+\sqrt{5}}{2},\quad \hat\phi=\frac{1-\sqrt{5}}{2}=1-\phi.
$$
Using the number $q=\hat\phi/\phi$, satisfying $-1<q<0$ and already defined in (\ref{eq:golden}),
leads to
\begin{equation}\label{eq:help}
F_n=\frac{1}{\sqrt{5}}\phi^n(1-q^n),\quad q\phi^2=-1,
\end{equation}
and for $\a\in\mathbb N$ and $n\ge 0$
$$
\frac{F_{\a}}{F_{\a+n}}=\frac{\sqrt{5}F_{\a}}{{\phi}^{\a+n}}\frac{1}{1-q^{\a+n}}
=(1-q^{\a})\sum_{k=0}^\infty (q^k/\phi)^n q^{\a k},
$$
which is the $n$'th moment of the real-valued  measure 
\begin{equation}\label{eq:mup}
\mu_{\a}=(1-q^{\a})\sum_{k=0}^\infty q^{\a k}\delta_{q^k/\phi}
\end{equation}
with total mass 1. When $\a$ is even then $\mu_{\a}$ is a probability
measure, but
when $\a$ is odd the masses $q^{\a k}$ change sign with the parity of $k$. Note
that $\mu_2$ is the measure considered in (\ref{eq:fibmea}).

For the Fibonomial coefficients defined in (\ref{eq:fibonomial}) one
has 
$$
\tbinom{n}{k}_{\mathbb F}=1,\; 0\le k\le n\le 2,
$$
and they satisfy a recursion formula 
\begin{equation}\label{eq:fibonomialrec}
\tbinom{n}{k}_{\mathbb F}=F_{k-1}\tbinom{n-1}{k}_{\mathbb F}
+F_{n-k+1}\tbinom{n-1}{k-1}_{\mathbb F},\;n>k\ge 1,
\end{equation}
see \cite{Kn}, which shows that the Fibonomial coefficients are integers.
From (\ref{eq:fibonomial}) it is also clear that
$$
\tbinom{n}{k}_{\mathbb F}=\tbinom{n}{n-k}_{\mathbb F},\quad 0\le k\le n.
$$
In \cite[Section 7.3]{G:R} one finds a discussion of the little
$q$-Jacobi polynomials defined in (\ref{eq:lqJacobi}), and it is
proved that
\begin{equation}\label{eq:ortJacobi}
\sum_{k=0}^\infty p_n(q^k;a,b;q)p_m(q^k;a,b;q)\frac{(bq;q)_k}{(q;q)_k}
(aq)^k=\frac{\delta_{n,m}}{h_n(a,b;q)},
\end{equation}
where
\begin{equation}\label{eq:norm}
h_n(a,b;q)=\frac{(abq;q)_n(1-abq^{2n+1})(aq;q)_n(aq;q)_\infty}
{(q;q)_n(1-abq)(bq;q)_n(abq^2;q)_\infty}(aq)^{-n}.
\end{equation}
In \cite{G:R} it is assumed that $0<q,aq<1$, but the derivation shows
that it holds for $|q|<1,|a|\le 1, |b|\le 1$, in particular in the
case of interest here: $-1<q<0,a=q^{\a-1}, b=1$, in the case of which we get  
\begin{equation}\label{eq:ortJacobispec}
\sum_{k=0}^\infty p_n(q^k;q^{\a-1},1;q)p_m(q^k;q^{\a-1},1;q)
q^{\a k}=\delta_{n,m}\frac{q^{\a n}(q;q)_n^2}{(q^{\a};q)_n^2(1-q^{\a+2n})}.
\end{equation}
This shows that the polynomials 
$$
p_n(x\phi;q^{\a-1},1;q)
$$
 are orthogonal with respect to $\mu_{\a}$
and that
$$
\int p_n(x\phi;q^{\a-1},1;q)p_m(x\phi;q^{\a-1},1;q)\,d\mu_{\a}(x)=
 \delta_{n,m}\frac{(1-q^{\a})q^{\a n}(q;q)_n^2}{(q^{\a};q)_n^2(1-q^{\a+2n})}.
$$
To simplify this apply (\ref{eq:help}) to  get
$$
\frac{(1-q^{\a})q^{\a n}(q;q)_n^2}{(q^{\a};q)_n^2(1-q^{\a+2n})}=
(-1)^{\a n}\tfrac{F_{\a}}{F_{\a+2n}}\left(\prod_{j=0}^{n-1}\tfrac{F_{1+j}}
{F_{\a+j}}\right)^2=(-1)^{\a
n}\tfrac{F_{\a}}{F_{\a+2n}}\tbinom{\a+n-1}{n}_{\mathbb F}^{-2}.
$$ 

\begin{cbthm}\label{cbthm:integercoef} Let $\a\in\mathbb N$.
The polynomials $p_n^{(\a)}(x)$ defined by
\begin{equation}\label{eq:poldef}
p_n^{(\a)}(x)=\tbinom{\a+n-1}{n}_{\mathbb F}p_n(x\phi;q^{\a-1},1;q)
\end{equation}
can be written
\begin{equation}\label{eq:pol}
p_n^{(\a)}(x)=\sum_{k=0}^n (-1)^{kn-\binom{k}{2}}\tbinom{n}{k}_{\mathbb F}
\tbinom{\a+n+k-1}{n}_{\mathbb F}x^k,
\end{equation}
and they satisfy
\begin{equation}\label{eq:pol1}
\int p_n^{(\a)}(x)p_m^{(\a)}(x)\,d\mu_{\a}(x)=\delta_{n,m}(-1)^{\a
  n}\tfrac{F_{\a}}{F_{\a+2n}},
\end{equation} 
so the corresponding orthonormal polynomials are
\begin{equation}\label{eq:polorth}
P_n^{(\a)}(x)=\sqrt{(-1)^{\a
  n}F_{\a+2n}/F_{\a}}p_n^{(\a)}(x).
\end{equation}
\end{cbthm}

{\it Proof}. By definition, see (\ref{eq:lqJacobi})
\begin{eqnarray*}
p_n^{(\a)}(x)&=&\tbinom{\a+n-1}{n}_{\mathbb F}\sum_{k=0}^n\frac{(q^{-n},q^{\a+n};q)_k}
{(q,q^{\a};q)_k}(q\phi
x)^k\\
&=&\tbinom{\a+n-1}{n}_{\mathbb F}\sum_{k=0}^n{\scriptsize
\left[\begin{matrix}n\\k\end{matrix}\right]_q}
\frac{(q^{\a+n};q)_k}{(q^{\a};q)_k}(-1)^kq^{\binom{k}{2}-nk}(q\phi x)^k,
\end{eqnarray*}
where 
$${\scriptsize
\left[\begin{matrix}n\\k\end{matrix}\right]_q}=
\frac{(q;q)_n}{(q;q)_k(q;q)_{n-k}}
$$
is the $q$-binomial coefficient.
Using (\ref{eq:help}) leads to
$${\scriptsize
\left[\begin{matrix}n\\k\end{matrix}\right]_q}=\binom{n}{k}_{\mathbb F}
\phi^{k(k-n)},
$$
hence
$$
p_n^{(\a)}(x)=\tbinom{\a+n-1}{n}_{\mathbb F}\sum_{k=0}^n(-1)^k\tbinom{n}{k}_{\mathbb F}
(\phi^2q)^{\binom{k+1}{2}-nk}\prod_{j=0}^{k-1}\tfrac{F_{\a+n+j}}{F_{\a+j}}x^k,
$$
which by (\ref{eq:help}) can be reduced to (\ref{eq:pol}).
 \quad$\square$

\begin{cbrem} {\rm The polynomials $p_n^{(\a)}(x)$ for $\a=1$ and $\a=2$
    are the polynomials in (\ref{eq:lqJacobi1}) and in (\ref{eq:fibpol1})
    respectively.}
\end{cbrem}

\begin{cbcor}\label{eq:det} For $\a\in\mathbb N$ 
$$
\det(1/F_{\a+i+j})_0^n=\left((-1)^{\a\binom{n+1}{2}}F_{\a}\prod_{k=1}^n
  F_{\a+2k}\tbinom{\a+2k-1}{k}_{\mathbb
    F}^2\right)^{-1},
$$
which is the reciprocal of an integer.
\end{cbcor}

{\it Proof}. 
From the general theory it is  known that the leading coefficient of the
orthonormal polynomial $P_n^{(\a)}$ is $\sqrt{D_{n-1}/D_n}$, where
$$
D_n=\det(F_{\a}/F_{\a+i+j})_0^n.
$$
From (\ref{eq:pol}) and (\ref{eq:polorth}) we then get
$$
D_{n-1}/D_n=(-1)^{\a n}\tfrac{F_{\a+2n}}{F_{\a}}
\tbinom{\a+2n-1}{n}_{\mathbb F}^2,
$$
hence
$$
\frac{1}{D_n}=\prod_{k=1}^n
\frac{D_{k-1}}{D_{k}}=(-1)^{\a\binom{n+1}{2}}\tfrac{1}{F_{\a}^n}
\prod_{k=1}^n F_{\a+2k} \tbinom{\a+2k-1}{k}_{\mathbb F}^2
$$
and the formula follows. \quad$\square$

\begin{cbthm}\label{cbthm:integer1} The $i,j$'th entry of
 the inverse of the  matrix $(1/F_{\a+i+j})_0^n$ is given as
\begin{equation}\label{eq:integer1}
(-1)^{n(\a+i+j)-\binom{i}{2}-\binom{j}{2}} F_{\a+i+j}
\tbinom{\a+n+i}{n-j}_{\mathbb F}\tbinom{\a+n+j}{n-i}_{\mathbb F}
\tbinom{\a+i+j-1}{i}_{\mathbb F}\tbinom{\a+i+j-1}{j}_{\mathbb F}.
\end{equation}
\end{cbthm}

{\it Proof.} From Theorem \ref{cbthm:A} we get
$$
\left(\left(F_{\a}/F_{\a+i+j}\right)_0^n\right)^{-1}=
\left(a_{i,j}^{(n)}(\a)\right)_0^n,
$$
where $a_{i,j}^{(n)}(\a)$ is the coefficient to $x^{i}y^{j}$ in the
kernel polynomial $K_n(x,y)$ for the orthonormal polynomials $P^{(\a)}_n$.
Inserting the expressions (\ref{eq:pol}) and
 (\ref{eq:polorth}) in the kernel polynomial and changing the order of
 summation gives
$$
F_{\a} a^{(n)}_{i,j}(\a)=\sum_{k=\max(i,j)}^n C^{(\a)}(k;i,j),
$$
where we for $k\ge i,j$  have defined
\begin{equation}\label{eq:sum}
C^{(\a)}(k;i,j):=(-1)^{k(\a+i+j)-\binom{i}{2}-\binom{j}{2}}
F_{\a+2k}\tbinom{k}{i}_{\mathbb F}\tbinom{k}{j}_{\mathbb F}
\tbinom{\a+k+i-1}{k}_{\mathbb F}\tbinom{\a+k+j-1}{k}_{\mathbb F}.
\end{equation}
To prove that this expression can be summed to give
(\ref{eq:integer1}), we use induction in $n$. By symmetry we can
always assume $i\ge j$. The starting step $n=k=i\ge j$ is easy and is
left to the reader. For the induction step let $R^{(\a)}(n;i,j)$ denote
the expression (\ref{eq:integer1}). It has to be established that
$$ 
R^{(\a)}(n+1;i,j)-R^{(\a)}(n;i,j)=C^{(\a)}(n+1;i,j).
$$
The left-hand side of this expression can be written
$$
(-1)^{(n+1)(\a+i+j)-\binom{i}{2}-\binom{j}{2}}F_{\a+i+j}
\tbinom{\a+i+j-1}{i}_{\mathbb F}\tbinom{\a+i+j-1}{j}_{\mathbb F} T,
$$
where
$$
T=\tbinom{\a+n+1+i}{n+1-j}_{\mathbb
  F}\tbinom{\a+n+1+j}{n+1-i}_{\mathbb F}-
(-1)^{\a+i+j}\tbinom{\a+n+i}{n-j}_{\mathbb F}\tbinom{\a+n+j}{n-i}_{\mathbb F}
$$
$$=\frac{(F_{\a+n+i}\cdots F_{\a+i+j+1})(F_{\a+n+j}\cdots F_{\a+i+j+1})}
{(F_1\cdots F_{n+1-j})(F_1\cdots F_{n+1-i})}\,\cdot
$$
$$
\left[F_{\a+n+i+1}F_{\a+n+j+1}-(-1)^{\a+i+j}
F_{n+1-i}F_{n+1-j}\right].
$$
By Lemma \ref{cbthm:fiblemma} below (with $n$ replaced by $n+1$),
 the expression in brackets equals
$F_{\a+2n+2}F_{\a+i+j}$, and now it is easy to complete the proof.
\quad$\square$

\begin{cblem}\label{cbthm:fiblemma} For $n\ge i,j\ge 0$ and $\a\ge 0$ the
 following formula holds
\begin{equation}\label{eq:fiblemma}
F_{\a+2n}F_{\a+i+j}=F_{\a+n+i}F_{\a+n+j}-(-1)^{\a+i+j}F_{n-i}F_{n-j}.
\end{equation}
\end{cblem}

{\it Proof.} Using Binet's formula, the
right-hand side of (\ref{eq:fiblemma}) multiplied with 5 equals
$$
(\phi^{\a+n+i}-{\hat\phi}^{\a+n+i})(\phi^{\a+n+j}-{\hat\phi}^{\a+n+j}) 
-(-1)^{\a+i+j}(\phi^{n-i}-{\hat\phi}^{n-i})(\phi^{n-j}-{\hat\phi}^{n-j}).
$$
Using $\phi\hat{\phi}=-1$ one gets after some simplification
$$ 
(\phi^{\a+2n}-{\hat\phi}^{\a+2n})(\phi^{\a+i+j}-{\hat\phi}^{\a+i+j}),
$$
which  establishes the formula. \quad$\square$

\begin{cbrem} {\rm For $\a=1$ the expression (\ref{eq:integer1}) reduces
to
$$
(-1)^{n(i+j+1)-\binom{i}{2}-\binom{j}{2}} F_{i+j+1}
\tbinom{n+i+1}{n-j}_{\mathbb F}\tbinom{n+j+1}{n-i}_{\mathbb F}
\tbinom{i+j}{i}_{\mathbb F}^2,
$$
which is the expression found by Richardson \cite{Ri}, except that he expressed
the sign in a different but equivalent manner.}
\end{cbrem}
\section{The Hilbert matrices}

For $\a>0$ the  matrices
\begin{equation}\label{eq:Hilbert}
\mathcal H_n^{(\a)}=\left(\a/(\a+i+j)\right)_0^n,\quad n=0,1,\ldots,
\end{equation}
are the Hankel matrices for the moment sequence
$$
s_n^{(\a)}=\a\int_0^1 x^nx^{\a-1}\,dx=\frac{\a}{\a+n},\quad n=0,1,\ldots
$$
of the measure $\sigma_{\a}=\a x^{\a-1}1_{]0,1[}(x)\,dx$.
The corresponding orthogonal polynomials are easily seen to be
\begin{equation}\label{eq:Legendre}
r_n^{(\a)}(x)=\frac{1}{n!}x^{-\alpha+1}D^n\;[x^{\a-1+n}(1-x)^n]=
(-1)^n\sum_{k=0}^n\tbinom{n}{k}\tbinom{\a-1+n}{k}(x-1)^k x^{n-k},
\end{equation}
since they are Jacobi polynomials transfered to $]0,1[$, cf. \cite{A:A:R}.
Using the binomial formula for $(x-1)^k$ we find
$$
r_n^{(\a)}(x)=(-1)^n\sum_{j=0}^n (-1)^j x^{n-j}c_j,
$$
where
\begin{eqnarray*}
c_j&=&\sum_{k=j}^n \tbinom{k}{j}\tbinom{n}{k}\tbinom{\a-1+n}{k}
=\sum_{l=0}^{n-j}\tbinom{j+l}{j}\tbinom{n}{j+l}\tbinom{\a-1+n}{j+l}\\
&=&\tbinom{n}{j}\tbinom{\a-1+n}{j}
 {}_2F_1{\scriptsize\left(\begin{matrix} -n+j,-n-\a+j+1\\
j+1\end{matrix};1\right)}=\tbinom{n}{j}\tbinom{2n+\a-j-1}{n},
\end{eqnarray*}
where the ${}_2F_1$ is summed by the Chu-Vandermonde formula,
cf. \cite[p. 67]{A:A:R}.
This gives
\begin{equation}\label{eq:Legendre1}
r_n^{(\a)}(x)=\sum_{j=0}^n(-1)^j\tbinom{n}{j}\tbinom{\a+n+j-1}{n}x^j.
\end{equation}
The orthonormal polynomials with positive leading coefficients
are given as
$$
R_n^{(\a)}(x)=(-1)^n\sqrt{\frac{\a+2n}{\a}}r_n^{(\a)}(x),
$$
so the corresponding kernel polynomials have  coefficients
$a_{i,j}^{(n)}(\a)$ which by Theorem \ref{cbthm:A} satisfy
\begin{equation}\label{eq:nyformel} 
\a a^{(n)}_{i,j}(\a)=(-1)^{i+j}\sum_{k=\max{(i,j)}}^n
(\a+2k)\tbinom{k}{i}\tbinom{k}{j}\tbinom{\a+k+i-1}{k}
\tbinom{\a+k+j-1}{k}.
\end{equation} 

\begin{cbthm}\label{cbthm:inversehilb} The $i,j$'th element of the inverse
  matrix of $\left(1/(\a+i+j)\right)_0^n$ is given as
\begin{equation}\label{eq:choi} 
(-1)^{i+j}(\a+i+j)\tbinom{\a+n+i}{n-j}\tbinom{\a+n+j}{n-i}
\tbinom{\a+i+j-1}{i}\tbinom{\a+i+j-1}{j}.
\end{equation}
In particular they are integers for $\a\in\mathbb N$. Furthermore,
\begin{equation}\label{eq:hilbalphadet}
\det\left(1/(\a+i+j)\right)_0^n=\left(\a\prod_{k=1}^n (\a+2k)
\tbinom{\a+2k-1}{k}^2\right)^{-1}.
\end{equation}
\end{cbthm}

{\it Proof}. Let $R(n;i,j)$ denote the number given in
(\ref{eq:choi}), and define 
$$
C(k;i,j)=(-1)^{i+j}(\a+2k)\tbinom{k}{i}\tbinom{k}{j}\tbinom{\a+k+i-1}{k}
\tbinom{\a+k+j-1}{k}, \quad k\ge i,j.
$$
We shall prove that
$$
R(n;i,j)=\sum_{k=\max(i,j)}^n C(k;i,j)
$$
by induction in $n$ and can assume $i\ge j$. This is easy for $n=k=i$
and we shall establish
\begin{equation}\label{eq:induc}
R(n+1;i,j)-R(n;i,j)=C(n+1;i,j).
\end{equation}
The left-hand side of this expression can  be written
$$
(-1)^{i+j}(\a+i+j)\tbinom{\a+i+j-1}{i}\tbinom{\a+i+j-1}{j}T,
$$
where
$$
T=\tbinom{\a+n+1+i}{n+1-j}\tbinom{\a+n+1+j}{n+1-i}-\tbinom{\a+n+i}{n-j}
\tbinom{\a+n+j}{n-i}
$$
$$
=\tfrac{((\a+n+i)\cdots(\a+i+j+1))((\a+n+j)\cdots(\a+i+j+1))}
{(n+1-j)!(n+1-i)!}\,\cdot
$$
$$
[(\a+n+1+i)(\a+n+1+j)-(n+1-j)(n+1-i)].
$$
The quantity in brackets equals $(\a+2n+2)(\a+i+j)$, and now it is
easy to complete the proof of (\ref{eq:induc}).

The leading coefficient of $R_n^{(\a)}(x)$ is
$$
\sqrt{\frac{D_{n-1}}{D_n}}=\sqrt{\frac{\a+2n}{\a}}\binom{\a+2n-1}{n},
$$
where 
$$
D_n=\det\left(\a/(\a+i+j)\right)_0^n=\a^{n+1}\det\left(1/(\a+i+j)\right)_0^n.
$$
Therefore
$$
\frac{1}{D_n}=\prod_{k=1}^n \frac{D_{k-1}}{D_k}=\frac{1}{\a^n}
\prod_{k=1}^n (\a+2k)\tbinom{\a+2k-1}{k}^2,
$$
which proves (\ref{eq:hilbalphadet}). \quad$\square$

\medskip

Replacing $x$ by $1-x$, we see that $r_n^{(\a)}(1-x)$ are 
orthogonal polynomials with respect to the probability measure
$\a(1-x)^{\a-1}1_{]0,1[}(x)\,dx$. The corresponding moment sequence is
\begin{equation}\label{eq:binomalpha}
s_n=\frac{1}{\binom{\a+n}{n}},
\end{equation}
and the corresponding orthonormal polynomials are 
$\sqrt{(\a+2n)/\a}\;r_n^{(\a)}(1-x)$. Therefore
\begin{equation}\label{eq:Knfinal}
 K_n(x,y)=\sum_{k=0}^n \frac{\a+2k}{\a}r_k^{(\a)}(1-x)r_k^{(\a)}(1-y),
\end{equation}
showing that the  coefficient to $x^{i}y^{j}$ in $\a K_n(x,y)$ is
an integer when $\a\in\mathbb N$. This yields

\begin{cbthm}\label{cbthm:binomalpha} Let $\a\in\mathbb N$. The inverse of
  the matrix
\begin{equation}\label{eq:hankelbinom}
\left(\frac{1}{\a\binom{\a+i+j}{\a}}\right)_0^n
\end{equation}
has integer entries.
\end{cbthm}

It is not difficult to prove that
$$
r_n^{(\a)}(1-x)=\sum_{k=0}^n
(-1)^{n-k}\tbinom{n}{k}\tbinom{\a+n+k-1}{k}x^k,
$$
and it follows that the entries of the inverse of
(\ref{eq:hankelbinom}) are given as
$$
(-1)^{i+j}\sum_{k=\max(i,j)}^n (\a+2k)\tbinom{k}{i}\tbinom{k}{j}
\tbinom{\a+k+i-1}{i}\tbinom{\a+k+j-1}{j}.
$$
This formula holds of course for any $\a>0$.

The results of this section for $\a=1,2$ have been treated in the survey paper
\cite{Be}, written in Danish. For $\a=1$ the formula for the
elements of the inverse  of $\mathcal H_n^{(\a)}$  was
given in \cite{Ch}, but goes at least back to Collar \cite{Co},
 while the formula for its determinant goes back to
Hilbert in \cite{Hi}. In this case the polynomials $r_n^{(1)}(x)$ are
the Legendre polynomials for the interval $[0,1]$,
cf. \cite[Section 7.7]{A:A:R}. These polynomials have succesfully been
used in the proof of the irrationality of $\zeta(3)$.
For $\a=2$ we have $(\a+2k)/\a=1+k$, so the coefficient to
$x^{i}y^{j}$ in (\ref{eq:Knfinal}) is an integer. In this case Theorem
\ref{cbthm:binomalpha} can be sharpened: The inverse of the matrix 
$\left(1/\tbinom{2+i+j}{2}\right)_0^n$ has integer coefficients. This
result is also given in \cite{Ri}.

\medskip
{\bf Added June 2007} A result equivalent to Theorem \ref{cbthm:A} is
given by Collar in \cite{Co}. Denoting by
$$
M_n=(p_{ij}),\quad 0\le i,j\le n
$$
the matrix of coefficients of the orthonormal polynomials, i.e.
$$
P_i(x)=\sum_{j=0}^n p_{ij}x^j,\quad i=0,1,\ldots,n,
$$
where $p_{ij}=0$ for $i<j$, then the orthonormality can be expressed as
the matrix equation  $M_nH_nM_n^t=I_n$, hence
\begin{equation}\label{eq:collar}
H_n^{-1}=M_n^t M_n.
\end{equation}
Collar uses (\ref{eq:collar}) to obtain formula (\ref{eq:choi})
 and states: \lq\lq Equation (\ref{eq:collar}), which provides an
 elegant method for the computation of the reciprocal of a moment
 matrix, is due to Dr A.\ C.\ Aitken. The author is grateful to
 Dr. Aitken for permission to describe the method and for many helpful
 suggestions.\rq\rq

The paper by Collar is not mentioned in Choi's paper \cite{Ch} and was
not included in the list of references in the first version of this paper.

In \cite{A:B} the authors have defined a $q$-analogue of the Hilbert
matrix for any complex $q$ different from the roots of unity and
have proved a $q$-analogue of (\ref{eq:choi}). When 
$q=(1-\sqrt{5})/(1+\sqrt{5})$ one can recover the results about
the Filbert matrices and for $q=-e^{-2\theta},\theta>0$  results of
Ismail (\cite{Is1}) about Hankel matrices of generalized Fibonacci
 numbers.

\vskip1truecm
\noindent Department of Mathematics\\
University of Copenhagen\\
Universitetsparken 5\\
DK-2100 K\o benhavn \O, Denmark

\medskip
\noindent Email: berg@math.ku.dk

 \end{document}